\documentclass[11 pt]{amsart}

\usepackage{amscd,amsmath,latexsym,amsbsy,amssymb,amscd}
\usepackage[lite]{amsrefs}
\usepackage{verbatim}
\usepackage{microtype}
\usepackage{mathrsfs}
\usepackage[mathscr]{euscript}
\usepackage[all]{xy}
\usepackage{layout}
\usepackage{enumerate}
\usepackage{dsfont}

\theoremstyle{plain}
\newtheorem{theorem}{Theorem}[section]
\newtheorem{lemma}[theorem]{Lemma}
\newtheorem{proposition}[theorem]{Proposition}

\newtheorem{corollary}[theorem]{Corollary}

\theoremstyle{definition}
\newtheorem{definition}[theorem]{Definition}

\theoremstyle{remark}

\DeclareMathOperator{\Hoh}{H}
\DeclareMathOperator{\Spin}{Spin}

\DeclareMathOperator{\Eoh}{E}
\DeclareMathOperator{\im}{im}

\DeclareMathOperator*{\colim}{colim}

\DeclareMathOperator{\Tor}{Tor}

\newcommand{\we}{\simeq}

\newcommand{\CC}{\mathds{C}}

\newcommand{\QQ}{\mathds{Q}}
\newcommand{\ZZ}{\mathds{Z}}
\newcommand{\PP}{\mathds{P}}
\newcommand{\Z}{\mathds{Z}}
\newcommand{\Mod}{\mathrm{Mod}}

\def\H{\mathcal{H}}
\def\J{\mathcal{J}}

\def\F{\mathcal{F}}

\def\limis{\displaystyle\mathop{\text{lim}^{1}}}
\def\lim{\displaystyle\mathop{\text{lim}}}
\def\holimi{\displaystyle\mathop{\text{holim}}}

\begin{document}

\title[Uniqueness of twisted K-theory]
{Actions of Eilenberg-MacLane spaces on K-theory spectra and uniqueness of twisted K-theory}

\author[B.~Antieau]{Benjamin Antieau$^{*}$}
\address{Department of Mathematics,
UCLA, 520 Portola Plaza, Los Angeles, CA 90095, USA}
\email{antieau@math.ucla.edu}
\thanks{$^{*}$The first author was supported in part 
by the NSF under Grant RTG
DMS 0838697}

\author[D.~Gepner]{David Gepner}
\address{Fakult\"at f\"ur Mathematik, 
Universit\"at Regensburg, 93040 Regensburg, Germany}
\email{djgepner@gmail.com}

\author[J. M. G\'omez]{Jos\'e Manuel G\'omez}
\address{Department of Mathematics,
Johns Hopkins University, Baltimore, MD 21218, USA}
\email{jgomez@math.jhu.edu}

\begin{abstract}
  \noindent
We prove the uniqueness of twisted $K$-theory in both the real 
and complex cases using the computation of the $K$-theories 
of Eilenberg-MacLane spaces due to Anderson and Hodgkin. 
As an application of our method, we give some vanishing results 
for actions of Eilenberg-MacLane spaces on $K$-theory spectra.
\end{abstract}

\date{\today}



\maketitle

\section{Introduction}

Twisted $K$-theory, due originally to Donovan and Karoubi 
\cite{donovan_karoubi}, has become an important concept bridging the
fields of analysis, geometry, topology and string theory. 
It is the home of many topological invariants which cannot be seen 
by untwisted $K$-theory in the same way that the fundamental class 
of a non-orientable manifold must live in twisted cohomology. 
For instance, an appropriate twisted $K$-theory is the receptor of a 
Thom-isomorphism for non-$\Spin^c$-vector bundles. Because of its 
importance and its place in many different fields, there are a wide 
variety of definitions appearing in the literature. In addition to 
\cite{donovan_karoubi}, there are the accounts of \cite{abg}, 
\cite{atiyah-segal}, \cite{atiyah_twisted_2006},  
\cite{bcmms}, \cite{fht},  \cite{karoubi}, \cite{may-sigurdsson},
\cite{rosenberg} and \cite{wang} to name just a few.
All of the these definitions exploit different models available 
for $K$-theory and it is a natural question to determine the 
relationship between all possible different approaches. 
The goal of this article is to show that all reasonable 
definitions of twisted (real or complex) $K$-theory 
essentially agree. 
 
In our context we take as a reasonable definition of
(real or complex) $K$-theory as one arising from a map 
$K(\ZZ,3)\rightarrow BGL_1 K$ or 
$K(\ZZ/2,2)\rightarrow BGL_1 KO$.
In general if $R$ is an $A_{\infty}$-spectrum  
we can twist the generalized cohomology represented by $R$ 
over a space $X$. Let $BGL_{1}R$ denote the classifying space 
of the space $GL_1 R$ of homotopy units in $R$. Given a 
map $f:X\to BGL_{1}R$ we have an induced map of 
$\infty$-categories 
$f_{*}:\Pi_\infty X\to \Mod_R$, where 
$\Pi_\infty $ is fundamental $\infty$-groupoid of $X$ and 
$\Mod_{R}$ is the  $\infty$-category of $R$-modules. We refer 
the reader to \cite{htt} for an account of $\infty$-categories.
This map factors through the full subgroupoid of $\Mod_{R}$ 
spanned by the free rank one $R$-module $R$. 
One can construct the $f$-twisted $R$-spectrum 
$R(X)_{f}$ by defining 
\[
R(X)_f=\colim_{\Pi_\infty X} f_{*}.
\]
This is an $R$-module that can be seen parametrized 
spectrum over $X$ whose fibers are rank one free $R$-modules. 
One can then define the $f$-twisted $R$-cohomology groups of $X$ 
by taking sections of this parametrized spectrum.
These ideas are made precise in \cite{abghr} and are outlined 
below in Section~\ref{sec:uniqueness}. 

We are particularly interested in the case where $R=K$ or $R=KO$,
the spectra representing real or complex $K$-theory. 
Thus, from the viewpoint of homotopy theory, there is only 
one definition of twisted $K$-theory: given a map 
$f:X\rightarrow BGL_1 K$ one produces the $f$-twisted  $K$-theory 
spectrum $K(X)_f$ over $X$. However, in applications, 
one typically wants to associate twists of $K$-theory arising 
from a geometrically accessible subspace of $BGL_1 K$. In the 
case of complex $K$-theory for example, we have an inclusion  
$K(\Z,3)\to BGL_1 K$ and we are interested in twists of 
$K$-theory arising from maps $X\to BGL_{1}K$ that factor 
through $K(\Z,3)$, at least up to homotopy. Such twists 
are classified by cohomology classes in $\Hoh^3(X,\ZZ)$. 
For instance, the definition of Donovan and Karoubi associates twisted 
$K$-theory spectra to the torsion classes in $\Hoh^3(X,\ZZ)$, 
while the definitions of Rosenberg and of Atiyah and Segal define 
twisted $K$-theory for all classes of $\Hoh^3(X,\ZZ)$.
In these cases a map $j:K(\ZZ,3)\rightarrow BGL_1 K$ is fixed, and
the $\alpha$-twisted $K$-theory, where $\alpha\in\Hoh^3(X,\ZZ)$ 
is a cohomology class classifying a map $f:X\to K(\ZZ,3)$, 
is the twisted $K$-theory corresponding to  the composition
\begin{equation*}
    X\xrightarrow{f}K(\ZZ,3)\xrightarrow{j}BGL_1 K.
\end{equation*}
The above construction depends on the map 
$j:K(\Z,3)\rightarrow BGL_1 K$ chosen and thus 
we face the  problem classifying maps 
$K(\ZZ,3)\rightarrow BGL_1 K$. 
Similarly, in the real case, one twists by elements of 
$\Hoh^2(X,\ZZ/2)$, and so desires a map 
$K(\ZZ/2)\rightarrow BGL_1 KO$. 
Our computations lead to the following theorem.

\begin{theorem} There are natural isomorphisms of groups
\begin{align*}
[K(\ZZ,3),BGL_1 K]\cong[K(\ZZ,3),K(\ZZ,3)]&\cong  \ZZ \\
[K(\ZZ/2,2),BGL_1 KO]\cong[K(\ZZ/2,2),K(\ZZ/2,2)]&\cong \ZZ/2.
\end{align*}
\end{theorem}

Thus, any two maps from $K(\ZZ,3)$ to $BGL_1 K$ differ by 
an endomorphism of $K(\ZZ,3)$, up to homotopy, and similarly 
any two maps from $K(\ZZ,2)$ to $BGL_1 KO$ differ by 
an endomorphism of $K(\ZZ/2,2)$, up to homotopy.

Next we outline our approach. By \cite{mst}, there is a 
decomposition of infinite loop spaces
\begin{equation*}
    GL_1 K\simeq K(\ZZ/2,0)\times K(\ZZ,2)\times BSU_{\otimes},
\end{equation*}
where $BSU_{\otimes}$ is the infinite loop space classifying 
virtual complex vector bundles of rank and determinant one, equipped 
with the tensor product structure. Since this splitting respects the 
infinite loop structures, it may be delooped, so that we obtain 
the splitting
\begin{equation*}
    BGL_1 K\simeq K(\ZZ/2,1)\times K(\ZZ,3)\times BBSU_{\otimes}.
\end{equation*}
Let $i:K(\ZZ,3)\rightarrow BGL_1 K$ be the canonical inclusion. 
We view a reasonable definition of $K$-theory as one arising 
from a map $j:K(\ZZ,3)\rightarrow BGL_1 K$, thus we wish to 
compare  $i$ and $j$.

Denote by $bsu_\otimes$ the connective  spectrum such that 
$\Omega^\infty bsu_\otimes\simeq BSU_{\otimes}$. 
The main result of this work says in the complex case that
\[
bsu^1_{\otimes}(K(\ZZ,3))=[K(\ZZ,3),BBSU_{\otimes}]=0.
\]
More generally,  in Section ~\ref{cohomology computations}
we provide conditions on a  finitely generated abelian group 
$\pi$ and $n$ that imply that  $bsu^1_{\otimes}(K(\pi,n))$ vanishes.
Our calculations rely on a computation of the $K$-theory of 
Eilenberg-MacLane spaces due to Anderson-Hodgkin 
\cite{anderson-hodgkin}. In particular, it follows that any 
map $K(\ZZ,3)\rightarrow BGL_1 K$ is homotopic to a integer 
multiple of $i$. In practice, to figure out which integer, it suffices 
to compute a differential in the twisted Atiyah-Hirzebruch 
spectral sequence, as is done in \cite{atiyah_twisted_2006}. 
All constructions appearing in the literature differ by a 
unit $\pm1$. 

In the real case,
\[
BGL_1 KO\simeq K(\ZZ/2,1)\times K(\ZZ/2,2)\times BBSO_{\otimes},
\]
and we show that
\[
bso^1_{\otimes}(K(\ZZ/2,2))=[K(\ZZ/2,2),BBSO_{\otimes}]=0.
\]
where $BSO_{\otimes}$ is the infinite loop space classifying 
virtual real vector bundles of rank and determinant one, equipped 
with the tensor product structure.  Here,  $bso_{\otimes}$ denotes the 
associated connective spectrum.

Therefore, any map $j:K(\ZZ/2,2)\rightarrow BGL_1 KO$ 
is either homotopically trivial or $j$ is equivalent to the canonical inclusion, and 
therefore there is a unique non-trivial definition of 
twisted real $K$-theory. This is proven in a similar way to the 
complex case.

One calls twists associated to a map from $X$ to $BBSU_\otimes$ 
(resp. to $BBSO_\otimes$) higher twists of $K$-theory on $X$. 
Thus our theorem amounts to saying that there are no higher twists 
of complex $K$-theory on $K(\ZZ,3)$ or of real $K$-theory on 
$K(\ZZ/2,2)$.  In fact, in Propositions \ref{prop:homtriv}, 
\ref{prop:homtriv real}, and~\ref{prop:sharpness} we determine
exactly when there are higher twists of $K$-theory on 
$K(\pi,n)$ for $\pi$ finitely generated and $n\geq 2$ 
(or $n\geq 3$ if $\pi$ is not torsion). The original result in this 
direction is due to the third named author \cite{gomez} who 
showed that there are no higher twists for complex $K$-theory 
on the classifying spaces of compact Lie groups $G$. 
The results in \cite{gomez} imply in particular that there 
are no higher twists of complex $K$-theory over $K(\pi,1)$ 
when $\pi$ is a finite group and also over $K(\Z^{n},2)$ for any 
$n\ge 0$ and thus our computations generalize these facts.

One might also be interested in twists of complex $K$-theory 
coming from $r$-torsion classes in $\Hoh^3(X,\ZZ)$ for some fixed 
integer $r$ as in \cite{antieau-williams}. We show
that $bsu_{\otimes}^1(K(\ZZ/r,2))=0$ so that the only twists of 
$K$-theory by $r$-torsion classes come from composing the 
Bockstein map  $\beta:K(\ZZ/r,2)\rightarrow K(\ZZ,3)$ with a 
map $K(\ZZ,3)\rightarrow BGL_1 K$.

The actions of Eilenberg-MacLane spectra appear as follows.
Given a map
\[
K(\pi,n)\rightarrow BGL_1 K,
\]
one may pass to the level of  loop spaces and obtain an $A_\infty$-map 
\[
K(\pi,n-1)\rightarrow GL_1 K. 
\]
For example, by looking at the map 
on loop spaces associated to  $i:K(\ZZ,3)\rightarrow BGL_1 K$
we obtain an $A_{\infty}$-map 
$\CC\PP^{\infty}\simeq K(\ZZ,2)\rightarrow GL_1 K$
which classifies the action of $\CC\PP^{\infty}$ on $K$ given by 
tensoring with line bundles. We call an $A_\infty$ map 
$K(\pi,n)\rightarrow GL_1 K$ an action of $K(\pi,n)$ on the 
$K$-theory spectrum. If the map factors through $BSU_{\otimes}$, 
we call the action a higher action. As a corollary of out computations 
we obtain  (Corollary~\ref{cor:actions}) the
classification of those $K(\pi,n)$ with finitely generated abelian 
group $\pi$ and $n\geq 2$ (or $n\geq 3$ if $\pi$ is not torsion) for which all higher actions are trivial.

\medskip

Here is the outline of the paper. The technical engine of the paper
is contained in Section~\ref{cohomology computations} where 
we compute the generalized cohomology groups $bsu^1_{\otimes}$ 
and $bso^1_{\otimes}$ of Eilenberg-MacLane spaces. After this,  
Section~\ref{sec:uniqueness} recalls the definition of
twisted $K$-theory via the method of \cite{abghr}, and we give 
there the proof of the uniqueness theorem.
Finally, in the appendix, we give a nice geometric model of 
$K$-theory in the spirit of Atiyah and Segal which has the 
advantage that it is a structured ring spectrum and so may be 
easily used to produce a map $K(\ZZ,3)\rightarrow BGL_1 K$.

\medskip

\paragraph{\bf{Notation:}} We will denote by $k$ the spectrum 
representing connective complex $K$-theory, by $K$ the spectrum 
representing complex $K$-theory, and by $KO$ the spectrum 
representing real $K$-theory. For a prime $p$ we will denote by 
$\Z_{p}$ the ring of $p$-adic integers. Given a spectrum $F$ and 
an abelian group $G$ we can introduce $G$ coefficients on $F$ 
by considering the spectrum $F_{G}=F\wedge MG$, where $MG$ 
is a Moore spectrum for the group $G$. Also, given an integer 
$n$, we denote the $(n-1)$-connected cover of $F$ by 
$F\left<n\right>$.\\

\paragraph{\bf{Acknowledgments:}} We would like to thank Ulrich Bunke 
for suggesting some very useful remarks regarding this problem and 
for making very detailed comments on an early draft. Also, we thank Peter 
Bousfield for some comments on the $K$-theory of 
Eilenberg-MacLane spaces.

\section{Cohomology computations}\label{cohomology computations}

The goal of this section is to determine when the groups 
$bsu_{\otimes}^1(K(\pi,n))$ and  $bso_{\otimes}^1(K(\pi,n))$ vanish. 
We show in detail the computations for the complex case.  
The real case is handled in a similar way and we only provide the main 
points leaving the details to the reader. 

\begin{lemma}\label{lem:ku-vanishing}
    Suppose that $\pi$ is a torsion abelian group. Then,
    \begin{equation*}
        \tilde{K}^{*}(K(\pi,n))=0
    \end{equation*}
    if $n\geq 2$. Suppose that $n\geq 3$ and that $\pi$ is non-torsion 
    (not necessarily torsion-free). Then,
    \begin{equation*}
        K^1(K(\pi,n))=0
    \end{equation*}
    if and only if $\pi\otimes_{\ZZ}\QQ$ is a $1$-dimensional 
    $\QQ$-vector space and $n$ is odd.
    \begin{proof}
        The case of $K(\pi,n)$ where $\pi$ is torsion and $n\geq 2$ is 
        \cite[Theorem 3]{yosimura}.
        Now, suppose that $\pi$ is non-torsion. Then,
        by \cite[Theorem 3]{yosimura},
        \begin{equation*}
            K^1(K(\pi,n))=K^1(K(\pi\otimes_{\ZZ}\QQ,n))=
            \bigoplus_{p+q=1}\Hoh^p\left(K(\pi\otimes_\ZZ\QQ,n),
            K^q(*)\right).
        \end{equation*}
        Therefore, it suffices to prove that if $\pi$ is a free abelian group
        and $n\geq 3$, then $K(\pi\otimes_{\ZZ}\QQ,n)$
        has integral cohomology concentrated in even degrees if and only if 
        $\pi\otimes_{\ZZ}\QQ$ is $1$-dimensional
        and $n$ is odd. We may as well assume that $\pi=\ZZ^I$ for some 
        non-empty set $I$. Let $\{e_i\}_{i\in I}$ be a basis for $\ZZ^I$.
        Since homology commutes with direct limits,
        by the computation of Cartan \cite[Th\'eor\`eme 1]{cartan-determination} of the 
        integral homology of Eilenberg-MacLane spaces,
        if $n$ is odd, then
        \begin{equation*}
            \Hoh_*(K(\ZZ^I,n),\QQ)\cong\Lambda_{\QQ}[\sigma^n e_i],
        \end{equation*}
        the rational exterior algebra on symbols $\sigma^n e_i$ in degree 
        $n$, and if $n$ is even, then
        \begin{equation*}
            \Hoh_*(K(\ZZ^I,n),\QQ)\cong \QQ[\sigma^n e_i],
        \end{equation*}
        a polynomial algebra. Since the reduced homology groups of 
        $K(\QQ^I,n)$ are $\QQ$-vector spaces, by universal 
        coefficients for homology,
        \begin{equation*}
            \tilde{\Hoh}_*(K(\QQ^I,n),\ZZ)\cong \tilde{\Hoh}_*(K(\QQ^I,n),\QQ).
        \end{equation*}
        On the other hand, $K(\ZZ^I,n)\rightarrow K(\QQ^I,n)$ is a 
        rational homotopy equivalence 
        (for instance, by \cite[Corollary~7.6]{griffiths-morgan}), so
        \begin{equation*}
            \tilde{\Hoh}_*(K(\ZZ^I,n),\QQ)\cong \tilde{\Hoh}_*(K(\QQ^I,n),\QQ).
        \end{equation*}
        Therefore, the reduced integral homology of $K(\QQ^I,n)$ is 
        concentrated in odd degrees if and only if $n$ is odd and $|I|=1$.
        But, since the reduced integral cohomology of $K(\QQ^I,n)$ 
        consists of $\QQ$-vector spaces, this implies that the reduced 
        integral cohomology of $K(\QQ^I,n)$ is concentrated in even 
        degrees if and only if $n$ is odd and $|I|=1$, by the universal 
        coefficient theorem.
    \end{proof}
\end{lemma}

\begin{lemma}\label{lem:ko-vanishing}
    Suppose that $\pi$ is a torsion abelian group. Then,
    \begin{equation*}
        \widetilde{KO}^*(K(\pi,n))=0
    \end{equation*}
    if $n\geq 2$.
    Suppose that $n\geq 3$ and that $\pi$ is non-torsion. Then,
    \begin{equation*}
        \widetilde{KO}^1(K(\pi,n))=0
    \end{equation*}
    if and only if $\pi\otimes_{\ZZ}\QQ$ is at most $3$-dimensional as 
    a $\QQ$-vector space and $n$ is odd.
    \begin{proof}
        If $\pi$ is torsion, then by the previous lemma,
        \begin{equation*}
            \widetilde{K}^*(K(\pi,n))=0
        \end{equation*}
        so that by \cite[Appendix]{anderson-hodgkin},
        \begin{equation*}
            \widetilde{KO}^*(K(\pi,n))=0.
        \end{equation*}
        In general,
        \begin{equation*}
            K(\pi,n)\rightarrow K(\pi\otimes_\ZZ\QQ,n)
        \end{equation*}
        induces an isomorphism on $K$-homology by \cite{yosimura},
        and therefore also an isomorphism on $KO$-cohomology 
        by \cite[Corollary~1.13]{meier}.
        Therefore, we can assume that $\pi$ is a free abelian group.
        Then, it is evidently sufficient to prove the statement for $\pi$ a 
        finitely generated free abelian group. Indeed if 
        $\widetilde{KO}^1(K(\tau,n))\neq 0$ for $\tau$ of rank at 
        least $4$, then choosing a splitting $\pi\cong \tau\oplus\sigma$ 
        shows that $\widetilde{KO}^1(K(\pi,n))\neq 0$ as well.
        Thus, let $\tau$ be a finitely generated free abelian group. 
        We show that $\widetilde{KO}^1(K(\tau,n))=0$
        if and only if $n\geq 3$ is odd and the rank of $\tau$ is at 
        most $3$. As we are in the finitely generated case, 
        by \cite[Appendix]{anderson-hodgkin},
        \begin{equation*}
            \widetilde{KO}^1(K(\tau,n))\cong 
            \bigoplus_{p+q=1}\Hoh^p(K(\tau\otimes_{\ZZ}\QQ,n),KO^q(*)).
        \end{equation*}
        Since the reduced cohomology of $K(\tau\otimes_\ZZ\QQ,n)$ is a 
        $\QQ$-vector space, in the direct sum above,
        \begin{equation*}
            \Hoh^p(K(\tau\otimes_\ZZ\QQ,n),KO^q(*))
        \end{equation*}
        can only be non-zero for $q=0\mod 4$. Therefore, 
        $\widetilde{KO}^1(K(\tau,n))=0$ if and only if
        $K(\tau\otimes_\ZZ\QQ,n)$ has no integral cohomology in 
        degrees equal to $1\mod 4$. If $n$ is even, then 
        $K(\tau\otimes_\ZZ\QQ)$ has integral cohomology in 
        degrees $1\mod 4$. Thus, 
        $\widetilde{KO}^1(K(\tau\otimes_\ZZ\QQ,n))\neq 0$.  If $n$ is odd, 
        supposing that $\tau=\ZZ^m$, the Cartan calculation we saw in 
        the proof of the previous lemma says that 
        $K(\tau\otimes_\ZZ\QQ,n)$ has integral cohomology in degrees
        \begin{equation*}
            n+1,2n+1,\ldots, mn+1.
        \end{equation*}
        If $m\leq 3$, then these degrees are $n+1,2n+1,3n+1$. 
        Since $n$ is odd, none of these numbers are equal to 
        $1\mod 4$. If $m\geq 4$, $4n+1=1\mod 4$. 
        This completes the proof.
    \end{proof}
\end{lemma}

Now we turn to the vanishing of $bsu_{\otimes}^1(K(\pi,n))$.
Let $\Sigma^4 k\simeq K\left<4\right>$ denote the $3$-connected cover of 
$K$-theory. This is a connective spectrum with infinite loop space 
$BSU_{\oplus}\simeq\Omega^\infty\Sigma^4 k$, though here the infinite loop 
space structure is additive and does not agree with the multiplicative one 
on $BSU_\otimes$. The main result of \cite[Corollary 1.4]{Adams}, however, 
asserts that the infinite loop structures become equivalent after 
localization or completion at any prime $p$. This implies in particular 
that $\Sigma^{4} k\wedge M\Z_{p}\simeq bsu_{\otimes}\wedge M\Z_{p}$ for 
every prime $p$. We are going to use this fact to show the triviality of 
$bsu_{\otimes}^{1}(K(\pi,n))$ for various $\pi$ and $n$. For this we need 
the following lemma.

\begin{lemma}\label{K theory K(Z,3)}
If (1) $\pi$ is a finite abelian group and $n\geq 2$ or (2) $\pi$ is a 
finitely generated abelian group with 
$\dim_{\QQ}\pi\otimes_{\ZZ}\QQ=1$ and $n\geq 3$ is odd, then 
$k_{\Z_p}^5(K(\pi,n))=0$ for every prime $p$ and $k^5(K(\pi,n))=0$.
    \begin{proof}
        By Lemma~\ref{lem:ku-vanishing} $K^1(K(\pi,n))=0$ in cases 
        (1) and (2). Note also that
        \begin{equation}\label{cohomology}
            \tilde{\Hoh}^{r}(K(\pi,n),\Z)=0 \text{ for }  0\le r\le 2
        \end{equation}
        under the hypotheses.
        
        Let $K(\pi,n)$ be endowed with a CW-complex structure with
        $m$-skeleton $F_m$, a finite CW-complex.  Note that 
        $\tilde{\Hoh}^r(F_{m},\ZZ)=0$ for $0\le r\le 2$ and $m$ large 
        enough. Then, by \cite{milnor}, there are exact sequences
        \begin{align*}
            0\rightarrow \limis_{m\to\infty}K^4(F_m)\rightarrow& K^5(K(\pi,n))
            \rightarrow\lim_{m\to\infty}K^5(F_m)\rightarrow 0,\\
            0\rightarrow \limis_{m\to\infty}k^4(F_m)\rightarrow & k^5(K(\pi,n))
            \rightarrow\lim_{m\to\infty}k^5(F_m)\rightarrow 0.
        \end{align*}
        We will prove separately that $\limis_{m\to\infty}k^4(F_m)=0$ 
        and $\lim_{m\to\infty}k^5(F_m)=0$.
        
        Let's show first that $\limis_{m\to\infty} k^4(F_m)=0$. Since 
        $K^5(K(\pi,n))=0$ we have $\limis_{m\to \infty} K^{4}(F_{m})=0$ 
        and it is easy to see that this implies that 
        $\limis_{m\to \infty} \tilde{K}^{4}(F_{m})=0$. Fix $m$ large enough 
        and consider the Atiyah-Hirzebruch spectral sequences computing 
        $K^{*}(F_{m})$ and $k^{*}(F_{m})$ 
        \begin{align}\label{eq:ahss}
            \Eoh_{2}^{r,s}=\Hoh^{r}(F_{m},k^{s}(*))
            &\Longrightarrow k^{r+s}(F_{m}),\\
            \tilde{\Eoh}_{2}^{r,s}=\Hoh^{r}(F_{m},K^{s}(*))
            &\Longrightarrow   K^{r+s}(F_{m}).
        \end{align}
        Both of these spectral sequences converge strongly as $F_{m}$ is a 
        finite CW-complex. We have a map of spectra $k\to K$ inducing an 
        isomorphism on homotopy groups in non-negative degrees. 
        This provides a map of spectral sequences 
        $f^{r,s}_{*}:\Eoh^{r,s}_{*}\to \tilde{\Eoh}^{r,s}_{*}$ such 
        that $f^{r,s}_{2}$ is an isomorphism whenever $s\le 0$. Moreover, 
        since $\tilde{\Hoh}^{r}(F_{m},\Z)=0$ for $0\le r\le 2$ we have that 
        $f^{r,s}_{2}$ is an isomorphism whenever $r+s=4$ and $r>0$. Also 
        note that there are no differentials that kill elements in total 
        degree $4$ in the case of $K$ that fail to do so in the case of $k$. 
        This is because the only possible such differentials must have source 
        $\tilde{\Eoh}_{*}^{1,2}$ but this is trivial as 
        $\tilde{\Eoh}_{2}^{1,2}=\Hoh^{1}(F_{m},\Z)=0$. 
        This proves that $f^{r,s}_{*}$ induces an isomorphism 
        $f^{r,s}_{\infty}:\Eoh_{\infty}^{r,s}\to \tilde{\Eoh}_{\infty}^{r,s}$ 
        whenever $r+s=4$ and $r>0$. Also 
        $f^{0,4}_{\infty}:\Eoh_{\infty}^{0,4}=0
        \to \tilde{\Eoh}_{\infty}^{0,4}\cong\Z$ 
        since $\tilde{\Eoh}_{2}^{0,4}=\Hoh^{0}(F_{m},\Z)\cong\Z$ and any 
        differential with source $\tilde{\Eoh}_{*}^{0,4}\cong\Z$ is trivial 
        as one sees by comparing $\tilde{\Eoh}_{*}^{r,s}$ with the 
        Atiyah-Hirzebruch spectral sequence computing $K(*)$. This in turn 
        proves that the map of spectra $k\to K$ induces a short exact sequence 
        \[
        0\to k^{4}(F_{m})\to K^{4}(F_{m})\to \Z\to 0.
        \] 
        Note that in fact $k^{4}(F_{m})\subset \tilde{K}^{4}(F_{m})$. 
        We conclude that the map $k\to K$ induces an isomorphism 
        $k^{4}(F_{m})\cong \tilde{K}^{4}(F_{m})$. 
        Since $\limis_{m\to \infty} \tilde{K}^{4}(F_{m})=0$ we conclude that 
        \[
        \limis_{m\to \infty} k^{4}(F_{m})=0.
        \]
        Let's prove now that $\lim_{m\to \infty}k^{5}(F_{m})=0$. To prove 
        this compare again the spectral sequences~\eqref{eq:ahss} 
        $\Eoh^{r,s}_{*}$ and $\tilde{\Eoh}^{r,s}_{*}$ for $F_m$. In total 
        degree $5$ the map $f^{r,s}_{*}$ is such that $f^{r,s}_{2}$ is an 
        isomorphism whenever $s\le 0$. A similar argument as before shows 
        also in this case there are no differentials that kill elements in 
        total degree $5$ in the case of $K$ that fail to do so in the case 
        of $k$. Therefore $f_{\infty}^{r,s}:\Eoh^{r,s}_{*}\to 
        \tilde{\Eoh}^{r,s}_{*}$ is an isomorphism whenever $r+s=5$ and 
        $s\le 0$. Also note that $\Eoh^{r,s}_{\infty}=0$ whenever 
        $r+s=5$ and $s> 0$. These facts show that 
        the map of spectra $k\to K$ induces an injective map 
        $k^{5}(F_{m})\to K^{5}(F_{m})$ for $m$ large enough. Indeed, a 
        map of filtered abelian groups with finite decreasing filtrations
        is injective if the map on each slice is injective by an iterated 
        use of the snake lemma. Given the commutative diagram 
        \begin{equation}
            \begin{CD}
                k^{5}(F_{m+1})@>{i_{*}}>>
                k^{5}(F_{m})\\
                @VVV @VVV \\
                K^{5}(F_{m+1})@>{i_{*}}>>
                K^{5}(F_{m}).
            \end{CD}
        \end{equation}
        and the fact that $\lim_{m\to \infty}K^{5}(F_{m})=0$, it follows 
        that $\lim_{m\to \infty}k^{5}(F_{m})=0$ since $\lim$ is left-exact. 
        The fact that
        $k_{\Z_{p}}^{5}(K(\pi,n))=0$ is proved in the same way once we 
        know that $K^5_{\Z_p}(K(\pi,n))=0$.
        To see this note that 
         \[
	      K_{\ZZ_{p}}=\holimi_{k\to \infty} K_{\ZZ/(p^{k})},
	      \]
	with structure maps coming from  the maps $\ZZ/(p^{k+1})\to \ZZ/(p^{k})$.  
	Because of this, we have a short exact sequence
	\begin{equation}\label{ses1}
	0\to \limis_{k\to \infty}K_{\ZZ/(p^{k})}^{4}(K(\pi,n))\to 
	K_{\ZZ_{p}}^{5}(K(\pi,n))\to 
	\lim_{k\to \infty}K_{\ZZ/(p^{k})}^{5}(K(\pi,n))\to 0.
	\end{equation}
         On the one hand, by \cite[Proposition 6.6]{Adamsbook} we have 
         a short exact sequence 
	\begin{equation}\label{eq8}
  0\to K^{5}(K(\pi,n))\otimes_{\ZZ} \ZZ/(p^{k})\to 
  K_{\ZZ/(p^{k})}^{5}(K(\pi,n)) \to
	\text{Tor}^{\ZZ}_{1}(K^{6}(K(\pi,n)),\ZZ/(p^{k}))\to 0.
  \end{equation}
        Under the given hypothesis $K^{5}(K(\pi,n))=0$ by 
        Lemma \ref{lem:ku-vanishing}. By \cite[Theorem I]{anderson-hodgkin} 
        we have $\tilde{K}^{*}(K(\pi,n))=0$ when $\pi$ is as in (1) and by 
        \cite[Theorem II]{anderson-hodgkin} we have that 
        $\tilde{K}^{*}(K(\pi,n))=H^{**}(K(\pi\otimes_\ZZ \QQ,n),\ZZ)=0$ when 
        $\pi$ satisfies (2).  In particular $\tilde{K}^{6}(K(\pi,n))$ is vector 
        space over $\QQ$ in this case. In either case it follows that  
        $\text{Tor}^{\ZZ}_{1}(K^{6}(K(\pi,n)),\ZZ/(p^{k}))=0$  and we 
        conclude that $K_{\ZZ/(p^{k})}^{5}(K(\pi,n))=0$.  This proves that the 
        right hand side in the short exact sequence (\ref{ses1}) 
        vanishes. We are left to prove that 
        \[
        \limis_{k\to \infty}K_{\ZZ/(p^{k})}^{4}(K(\pi,n))=0.
        \] 
        To show this we use the exact sequence
        \begin{equation}\label{ses2}
        0\to K^{4}(K(\pi,n))\otimes \ZZ/(p^{k})\to K_{\ZZ/(p^{k})}^{4}(K(\pi,n))
        \to \text{Tor}^{\ZZ}_{1}(K^{5}(K(\pi,n)),\ZZ/(p^{k}))\to 0.
        \end{equation}
        Since $K^{5}(K(\pi,n))=0$, we conclude from (\ref{ses2}) that 
        \[
        K_{\ZZ/(p^{k})}^{4}(K(\pi,n))=K^{4}(K(\pi,n))\otimes_\ZZ \ZZ/(p^{k}).
        \] 
        From here we can see that the maps 
        $K_{\ZZ/(p^{k+1})}^{4}(K(\pi,n))\to K_{\ZZ/(p^{k})}^{4}(K(\pi,n))$ 
        are surjective and thus the $\text{lim}^{1}$ term in the short exact 
        sequence (\ref{ses1}) vanishes. This proves that 
        $K_{\Z_p}^5(K(\pi,n))=0$.
    \end{proof}
\end{lemma}

A similar computation can be done in the real case. Consider $KO\left<2\right>$, 
the $1$-connected cover of $KO$. Then $KO\left<2\right>$ is a connective 
spectrum with $\Omega^{\infty}KO\left<2\right>\simeq BSO_{\oplus}$. By 
\cite[Corollary 1.4]{Adams} it follows that  
$KO\left<2\right>\wedge M\Z_{p}\simeq bso_{\otimes}\wedge M\Z_{p}$ for 
every prime $p$. 

\begin{lemma}
If (1) $\pi$ is a finite abelian group and $n\geq 2$ or (2) $\pi$ is a 
finitely generated abelian group with 
$1\leq\dim_{\QQ}\pi\otimes_{\ZZ}\QQ\leq 3$ and 
$n\geq 3$ is odd, then $KO\left<2\right>_{\Z_{p}}^{1}(K(\pi,n))=0$ 
for every prime $p$ and $KO\left< 2\right>^{1}(K(\pi,n))=0$.
\begin{proof}
        Let $\{F_{m}\}_{m\ge 0}$ be the skeleton filtration of CW-complex 
        structure on $K(\pi, n)$ in such a way that $F_{m}$ is a finite 
        CW-complex. Note that in these cases we also have for large $m$ 
        \begin{equation*}
            \tilde{\Hoh}^{r}(F_{m},\Z)=0 \text{ for }  0\le r\le 2.
        \end{equation*}
        Also, $\widetilde{KO}^{*}(K(\pi,n))=0$ 
        when $\pi$ is finite abelian and $n\geq 2$, and 
        $\widetilde{KO}^1(K(\pi,n))=0$ for $\pi$ finitely generated and 
        $n\geq 3$ odd, as proved above. We argue in a similar way 
        as in the previous lemma. We can compare the Atiyah-Hirzebruch 
        spectral sequences computing $KO\left<2\right>^{*}(F_{m})$ and 
        $KO^*(F_{m})$. By doing so we prove that 
        \[
        \limis_{m\to \infty}KO\left<2\right>^{0}(F_{m})=0 \text{ and } 
        \lim_{m\to \infty}KO\left<2\right>^{1}(F_{m})=0. 
        \]
        The lemma follows using the $\limis$ exact sequence in 
        $KO\left<2\right>^{1}$ associated to the filtration 
        $\{F_{m}\}_{m\ge 0}$.  The argument for $p$-completed $KO$-theory 
        of $K(\pi,n)$ follows the same lines as the complex case 
        using the fact that 
        \[
	      KO\left<2\right>_{\ZZ_{p}}=\holimi_{k\to \infty} 
	      KO\left<2\right>_{\ZZ/(p^{k})}.
        \]			 
	\end{proof}
\end{lemma}

\begin{definition}An inverse system of groups $\{G_n\}$, i.e., a diagram 
of the form
    \[
    \cdots\to G_{n+1}\to G_n\to\cdots \to G_{2}\to G_{1},
    \] 
    is said to satisfy the Mittag-Leffler condition if 
    for every $i$ we can find a $j>i$ such that for every $k>j$, 
    \[
    \im(G_{k}\to G_{i})=\im(G_{j}\to G_{i}).
    \]
\end{definition}

It is well known \cite[Proposition~3.5.7]{weibel} that if $\{G_{n}\}$ 
satisfies the Mittag-Leffler condition then 
\[
\limis_{k\to \infty} G_{k}=0.
\] 
On the other hand, if each $G_{k}$ is a countable group and 
$\limis_{m\to \infty} G_{m}=0$, then by \cite[Theorem 2]{McGibbon} we 
have that the system $\{G_{n}\}$ must satisfy the Mittag-Leffler 
condition.

\medskip

Using the previous propositions we can 
show that  $bsu_{\otimes}^{1}(K(\pi,n))$ vanishes for certain 
vales of $n$ and some abelian groups $\pi$. This computation is 
central for our treatment on twisted $K$-theory. 

\begin{proposition}\label{prop:homtriv}
    If (1) $\pi$ is a finite abelian group and $n\geq 2$ or (2) $\pi$ is a 
    finitely generated abelian group with 
    $\dim_{\QQ}\pi\otimes_{\ZZ}\QQ=1$ and $n\geq 3$ is odd, 
    then we have
    \[
    bsu_{\otimes}^{1}(K(\pi,n))=0.
    \]
    In particular, there are no higher twists of complex 
    $K$-theory  on $K(\pi,n)$ in either case.
    \begin{proof}
        As above, since $\pi$ is finitely generated, we can give 
        $K(\pi,n)$ a CW-complex structure in such a way that 
        the $m$-skeleton $F_{m}$ is a finite CW-complex. 
        The filtration $\{F_{m}\}_{m\ge 0}$ induces a short exact sequence 
        \[
        0\to \limis_{m\to \infty} bsu_{\otimes}^{0}(F_{m})
        \to bsu_{\otimes}^{1}(K(\pi,n))\to 
        \lim_{m\to \infty}bsu_{\otimes}^{1}(F_{m})\to 0.
        \]
        Below we prove that the sequence of groups 
        $\{bsu_{\otimes}^{0}(F_{m})\}_{m\ge 0}$ satisfies the 
        Mittag-Leffler condition. Thus the $\limis$ in the 
        previous sequence vanishes yielding
        \begin{equation}\label{equation1}
        bsu_{\otimes}^{1}(K(\pi,n))
        =\lim_{m\to \infty}bsu_{\otimes}^{1}(F_{m}).
        \end{equation}
        We are going to show that 
        $\lim_{m\to \infty}(bsu_{\otimes}^{1}(F_{m})\otimes \Z_{p})=0$ 
        for all prime numbers $p$. Because each $bsu_{\otimes}^1(F_m)$ 
        is finitely generated group then by \cite[Lemma 7]{gomez} 
        $\lim_{m\to \infty}bsu_{\otimes}^{1}(F_{m})=0$ and thus 
        $bsu_{\otimes}^{1}(K(\pi,n))=0$ by \eqref{equation1}.

        Consider the short exact sequence 
        \begin{equation}\label{lim1 for K(Z,3)}
            0\to \limis_{m\to \infty}(bsu_{\otimes}
            \wedge M\Z_{p})^{0}(F_{m})
            \to (bsu_{\otimes}\wedge M\Z_{p})^{1}(K(\pi,n))
            \to\lim_{m\to \infty}(bsu_{\otimes}
            \wedge M\Z_{p})^{1}(F_{m})\to 0.
        \end{equation}
        By \cite[Corollary 1.4]{Adams} we have 
        $\Sigma^{4} k\wedge M\Z_{p}\simeq bsu_{\otimes}\wedge M\Z_{p}$. 
        This together with Lemma \ref{K theory K(Z,3)} gives
        \[
        (bsu_{\otimes}\wedge M\Z_{p})^{1}(K(\pi,n))\cong 
        k_{\Z_{p}}^{5}(K(\pi,n))=0. 
        \]
        We conclude that the middle term in (\ref{lim1 for K(Z,3)}) 
        vanishes and hence we see that 
        \[
        \lim_{m\to \infty}(bsu_{\otimes}\wedge M\Z_{p})^{1}(F_{m})=0.
        \]
        By \cite[Proposition III.6.6]{Adamsbook} there is a short exact 
        sequence 
        \[
        0\to bsu_{\otimes}^{1}(F_{m})\otimes \Z_{p}\to (bsu_{\otimes}
        \wedge M\Z_{p})^{1}(F_{m})
        \to \Tor^{1}_{\ZZ}(bsu_{\otimes}^{2}(F_{m}),\Z_{p})\to 0.
        \]
        The $\Tor$ term in this sequence vanishes as $\Z_{p}$ is flat as 
        a $\Z$-module. Therefore 
        \[
        (bsu_{\otimes}\wedge M\Z_{p})^{1}(F_{m})
        =bsu_{\otimes}^{1}(F_{m})\otimes_\ZZ \Z_{p}
        \]  
        and this in turn shows that for every prime $p$
        \[
        \lim_{m\to \infty}(bsu_{\otimes}^{1}(F_{m})\otimes_\ZZ \Z_{p})=0.
        \]
        We are left to prove that the system of groups 
        $B_{m}:=bsu_{\otimes}^{0}(F_{m})$ satisfies the 
        Mittag-Leffler condition. For every $m\ge 0$ let $A_{m}:=k^{4}(F_{m})$. 
        We saw in the proof of Lemma~\ref{K theory K(Z,3)} that
        \begin{equation*}
            \limis_{m\to\infty}A_m=0.
        \end{equation*}
        Also since $F_{m}$ is a finite CW-complex, we have that 
        $A_{m}$ and $B_{m}$ are finitely generated abelian groups 
        for every $m\ge 0$, in particular they are countable. 
        Therefore the system $\{A_{m}\}_{m\ge 0}$ satisfies 
        the Mittag-Leffler condition.

        On the other hand, as pointed out above 
        $\Sigma^{4} k\wedge M\Z_{p}\simeq bsu_{\otimes}\wedge M\Z_{p}$, thus 
        for every $m\ge 0$ and any prime number $p$
        \begin{equation}\label{eq17}
        A_{m}\otimes_\ZZ \Z_{p}=k_{\Z_{p}}^{4}(F_{m})
        \stackrel{\simeq}{\rightarrow}
        (bsu_{\otimes}\wedge M\Z_{p})^{0}(F_{m})=B_{m}\otimes \Z_{p}.
        \end{equation}
        The outer equalities follow by \cite[Proposition III.6.6]{Adamsbook} 
        since the $\Tor$ terms also vanish here. This yields a commutative 
        diagram in which the vertical arrows are isomorphisms
        \begin{equation*}
        \begin{array}{ccccccc}
        \rightarrow  & A_{m}\otimes_\ZZ \Z_{p} & \cdots  & \rightarrow  
        & A_{2}\otimes_\ZZ
        \Z_{p} & \rightarrow  & A_{1}\otimes_\ZZ \Z_{p} \\ 
        & \downarrow  &  &  & \downarrow  &  & \downarrow  \\ 
        \rightarrow  & B_{m}\otimes_\ZZ \Z_{p} & \cdots  & \rightarrow  
        & B_{2}\otimes_\ZZ
        \Z_{p} & \rightarrow  & B_{1}\otimes_\ZZ \Z_{p}.%
        \end{array}%
        \end{equation*}%
        Using this diagram and the fact that $\{A_{m}\}_{m\ge 0}$ satisfies 
        the Mittag-Leffler property it can be seen that the system 
        $\{B_{m}\}_{m\ge 0}$ also satisfies the Mittag-Leffler condition 
        using an argument similar to that in \cite[Theorem 5]{gomez}.
    \end{proof}
\end{proposition}

The previous proposition has the following real analogue that can be 
proved in the same way using the fact that 
$KO\left<2\right>\wedge M\Z_{p}\simeq bso_{\otimes}\wedge M\Z_{p}$ for 
every prime $p$. 

\begin{proposition}\label{prop:homtriv real}
    If (1) $\pi$ is a finite abelian group and $n\geq 2$ or (2) $\pi$ is a 
    finitely generated abelian group
    with $1\leq\dim_{\QQ}\pi\otimes_{\ZZ}\QQ\leq 3$ and $n\geq 3$ is 
    odd, then we have
    \[
    bso_{\otimes}^{1}(K(\pi,n))=0.
    \]
    In particular, there are no higher twists of real 
    $K$-theory on $K(\pi,n)$ in either case.
\end{proposition}

On the other hand, Proposition~\ref{prop:homtriv} is sharp as 
we show next. A real analogue can 
be proved in a similar way. 

\begin{proposition}\label{prop:sharpness}
    If (1) $\pi$ is a non-torsion (not necessarily torsion-free) finitely 
    generated abelian group and $n>3$ is even or (2) $\pi$ is a 
    finitely generated abelian group with 
    $\dim_{\QQ}\pi\otimes_{\ZZ}\QQ>1$ and $n\geq 3$ is odd, then
    $bsu_{\otimes}^1(K(\pi,n))\neq 0$.
    \begin{proof}
        By Lemma~\ref{lem:ku-vanishing} we know that 
        $K^{5}(K(\pi,n))\ne 0$ in these cases.  Let's show first that 
        $k^{5}(K(\pi,n))\ne 0$. Assume by contradiction that 
        $k^{5}(K(\pi,n))= 0$. As before give  $K(\pi,n)$ a 
        structure of a CW-complex such that $F_{k}$, the $k$-skeleton 
        of $K(\pi,n)$, is a finite CW-complex. 
        Since we are assuming that $k^{5}(K(\pi,n))= 0$ we have that 
        $\limis_{k\to \infty}k^{4}(F_{k})=0$ and 
        $\lim_{k\to \infty}k^{5}(F_{k})=0$.
        By comparing the Atiyah-Hirzebruch spectral sequences computing 
        $K^{*}(F_k)$ and $k^{*}(F_k)$ as in Lemma \ref{K theory K(Z,3)} 
        we can see that $\limis_{k\to \infty}K^{4}(F_{k})=0$ and 
        $\lim_{k\to \infty}K^{5}(F_{k})=0$. This in turn proves that 
        $K^{5}(K(\pi,n))= 0$ which is a contradiction.
        Let's show now that $bsu_{\otimes}^{1}(K(\pi,n))\ne 0$. 
        Reasoning by contradiction again assume that 
        $bsu_{\otimes}^{1}(K(\pi,n))=0$. The short exact sequence 
        \[
        0\to \limis_{k\to \infty} bsu_{\otimes}^{0}(F_{k})
        \to bsu_{\otimes}^{1}(K(\pi,n))\to 
        \lim_{k\to \infty}bsu_{\otimes}^{1}(F_{k})\to 0
        \]
        shows that $\lim_{k\to \infty}bsu_{\otimes}^{1}(F_{k})=0$ 
        and $\limis_{k\to \infty} bsu_{\otimes}^{0}(F_{k})=0$. 
        Since $F_{k}$ is a finite CW-complex we have that 
        $B_{k}:=bsu_{\otimes}^{0}(F_{k})$ is finitely generated 
        for every $k\ge 0$. In particular we conclude that the 
        system $\{B_{k}\}_{k\ge 0}$ satisfies the Mittag-Leffler 
        condition.  Let $A_{k}=k^{4}(F_{k})$. As in the proof of the 
        previous proposition we have a commutative 
        diagram in which the vertical arrows are isomorphisms
        \begin{equation*}
        \begin{array}{ccccccc}
        \rightarrow  & A_{n}\otimes_\ZZ \Z_{p} & \cdots  & \rightarrow  
        & A_{2}\otimes_\ZZ
        \Z_{p} & \rightarrow  & A_{1}\otimes_\ZZ \Z_{p} \\ 
        & \downarrow  &  &  & \downarrow  &  & \downarrow  \\ 
        \rightarrow  & B_{n}\otimes_\ZZ \Z_{p} & \cdots  & \rightarrow  
        & B_{2}\otimes_\ZZ
        \Z_{p} & \rightarrow  & B_{1}\otimes_\ZZ \Z_{p}.%
        \end{array}%
        \end{equation*}%
        This diagram, the fact that $\{B_{k}\}_{k\ge 0}$ satisfies the 
        Mittag-Leffler property and an argument similar to the one 
        provided in \cite[Theorem 5]{gomez} prove that 
        $\{A_{k}\}_{k\ge 0}$ also satisfies the Mittag-Leffler 
        property and in particular 
        \[
        \limis_{k\to \infty}A_{k}=\limis_{k\to \infty}k^{4}(F_{k})=0.     
        \]  
        On the other hand, by \cite[Theorem III]{anderson-hodgkin} 
        we have $\lim_{k\to \infty}K^{5}(F_{k})=0$. Comparing 
        the Atiyah-Hirzebruch spectral sequences computing 
        $K^{*}(F_{k})$ and $k^{*}(F_{k})$ we can see that 
        $K^{5}(F_{k})\cong k^{5}(F_{k})$, in particular 
        we obtain $\lim_{k\to \infty}k^{5}(F_{k})=0$. Finally, the 
        short exact sequence 
        \[
        0\to \limis_{k\to \infty} k^{4}(F_{k})
        \to k^{5}(K(\pi,n))\to 
        \lim_{k\to \infty}k^{5}(F_{k})\to 0
        \]
        shows that $k^{5}(K(\pi,n))=0$ which is a contradiction.
        \end{proof}
\end{proposition}

Next we consider  actions of  Eilenberg-MacLane spaces on the 
$K$-theory spectrum. We call an $A_\infty$-map 
$K(\pi,n-1)\rightarrow GL_1 K$ an action of $K(\pi,n-1)$ on $K$.
Given given an $A_{\infty}$-action 
of $K(\pi,n-1)$ on $K$, it can be de-looped to obtain a map 
$K(\pi,n)\rightarrow BGL_1 K$. Conversely,  given a map 
$K(\pi,n)\rightarrow BGL_1 K$, we obtain an $A_\infty$-map 
$K(\pi,n-1)\rightarrow GL_1 K$ by passing to the level of loop space. 
In fact actions of $K(\pi,n-1)$ on the 
$K$-theory spectrum are in one to one correspondence with 
maps $K(\pi,n)\to  BGL_1 K$. We call an action of $K(\pi,n-1)$ on $K$
a higher action if the corresponding map $ K(\pi,n)\rightarrow BGL_1 K$ 
factors through $BBSU_{\otimes}$. The above work can be 
rephrased as follows.

\begin{corollary}\label{cor:actions}
    Let $\pi$ be a finitely generated abelian group and $n\geq 2$ an integer.
    \begin{enumerate}
        \item[\rm{(2)}] There are no higher actions of $K(\pi,n)$ on $K$
            if and only if $\pi$ is torsion or $n$ is even and
            $\dim\pi\otimes_\ZZ\QQ=1$.
        \item[\rm{(3)}] There are no higher actions of $K(\pi,n)$ on $KO$
            if and only if $\pi$ is torsion or $n$ is even and
            $1\leq\dim\pi\otimes_\ZZ\QQ\leq 3$.
    \end{enumerate}
\end{corollary}

\begin{corollary}\label{cor:finactions}
    Let $\pi$ be a finite abelian group. Then, there are no higher actions of $K(\pi,1)$ on
    either $K$ or $KO$.
\end{corollary}

Corollary~\ref{cor:finactions} was obtained by Gomez~\cite{gomez}.

\section{Uniqueness of twisted $K$-theory}\label{sec:uniqueness}

In this section we use the computations of $bsu_{\otimes}$ and 
$bso_{\otimes}$-cohomology in the previous section to establish a 
uniqueness statement for definitions of twisted $K$-theory for 
both the real and complex cases. 

Let $R$ denote an $A_{\infty}$-ring spectrum. We can twist 
the generalized cohomology represented by $R$ over a space $X$. 
Let $GL_{1}R$ the space of homotopy units of $R$. This space 
is defined as the homotopy pullback in 
the diagram
\[
\begin{CD}
GL_{1}R@>>> \Omega^{\infty}R\\
@VVV @VVV \\
(\pi_{0}\Omega R)^{\times}@>>>\pi_{0}\Omega^{\infty}R.    
\end{CD}
\]
The space $GL_{1}R$ is a group-like $A_{\infty}$-space and 
twists of the theory $R$ over a space $X$ are classified by 
homotopy classes of maps $X\to BGL_{1}R$. Specifically, given a map 
$f:X\rightarrow BGL_1 R$, we obtain 
(by associating to a space $X$ its fundamental $\infty$-groupoid 
$\Pi_{\infty} X$, as in \cite{htt}) an induced map of 
$\infty$-categories
\begin{equation*}
    \Pi_\infty X\stackrel{f}{\longrightarrow}
    \Pi_\infty BGL_1 R\simeq B\mathrm{Aut}_R(R)
    \stackrel{i}{\longrightarrow}\Mod_R.
\end{equation*}
Here $\mathrm{Aut}_R(R)$ denotes the group-like $A_\infty$-space 
of automorphisms of $R$ in $\Mod_R$, $B\mathrm{Aut}_R(R)$ its 
delooping, and $B\mathrm{Aut}_R(R)\to\Mod_R$ the inclusion of the 
full subgroupoid of $\Mod_R$ spanned by the free rank one $R$-module 
$R$.

The $R$-module spectrum $R(X)_f$, or the $f$-twisted $R$-theory 
spectrum of $X$, is the resulting ``Thom spectrum''
\[
R(X)_f=\colim_{\Pi_\infty X} i\circ f,
\]
the colimit in $\Mod_R$ of the composite map 
$i\circ f:\Pi_\infty X\to\Mod_R$. The colimit exists since $\Mod_R$ 
admits colimits indexed by an arbitrary small $\infty$-category. 
See \cite{htt} for an account of the $\infty$-categorical 
theory of colimits. The $f$-twisted $R$-cohomology groups of $X$ 
are 
\[
R^{n}(X)_{f}:=\pi_{0}F_{R}(R(X)_{f},\Sigma^n R),
\]
where $F_{R}(R(X)_{f},\Sigma^n R)$ is the function spectrum of 
$R$-module maps $R(X)_{f}\to \Sigma^n R$.

We can use this method of twisting in the particular cases of 
twisted (real or complex) $K$-theory.
In the complex case the decomposition
\[
GL_1 K\simeq K(\ZZ/2,0)\times K(\ZZ,2)\times BSU_\otimes
\]
is compatible with the evident $A_\infty$-structures so it 
deloops to a decomposition
\[
BGL_1 K\simeq K(\ZZ/2,1)\times K(\ZZ,3)\times BBSU_\otimes.
\]
Therefore, twists of complex $K$-theory over $X$ are 
classified by homotopy classes of maps 
\[
X\to BGL_1 K\simeq K(\ZZ/2,1)\times K(\ZZ,3)\times BBSU_\otimes.
\]
The twisted cohomology groups $K^{n}(X)_{f}$ depend only on the 
homotopy class of $f$, through non-canonical isomorphisms, and 
thus by the decomposition
\begin{equation*}
    BGL_1 K\we K(\ZZ/2,1)\times K(\ZZ,3)\times BBSU_{\otimes}
\end{equation*}
we have twists of complex $K$-theory associated to elements in 
\[
\Hoh^1(X,\ZZ/2)\times \Hoh^3(X,\ZZ)\times bsu_{\otimes}^{1}(X).
\]

In the geometric applications however, one specializes to twists of 
$K$-theory associated to maps $f$ representing a cohomology class in 
$\Hoh^3(X,\ZZ)$. More precisely, let 
\[
i:K(\ZZ,3)\longrightarrow BGL_1 K.
\]
be the inclusion map and suppose that $f:X\to K(\ZZ,3)$ is a 
map representing a cohomology class in $\Hoh^3(X,\ZZ)$. 
Then the $f$-twisted $K$-theory spectrum of $X$ is defined as 
$K(X)_{i\circ f}$.

Different constructions or models of twisted $K$-theory 
associated to cohomology classes in $\Hoh^3(X,\ZZ)$ can be 
constructed by specifying a map 
\[
j:K(\ZZ,3)\longrightarrow BGL_1 K.
\]
Given such a map we can define the twisted $K$-groups as above. 

Thus a particular 
model or definition of twisted (complex) $K$-theory 
associated to cohomology classes in $\Hoh^3(X,\ZZ)$
amounts to producing a particular map 
$K(\ZZ,3)\longrightarrow BGL_1 K$.
We discuss here how the construction given in 
\cite{atiyah_twisted_2006} by Atiyah and Segal fits into this 
framework. Let $\H$ be a fixed infinite dimensional separable 
Hilbert space. The space of Fredholm operators $Fred(\H)$ with 
the norm topology is then a classifying space for complex $K$-theory. 
The space of unitary operators $U(\H)$ acts by conjugation on 
$Fred(\H)$. This induces an action of the projective unitary 
group $PU(\H)$ on $Fred(\H)$. The space $PU(\H)$ is a $K(\ZZ,2)$ 
and given a map 
\[
f:X\to BPU(\H)\simeq K(\ZZ,3)
\] 
there is an associated principal $PU(\H)$-bundle $P\to X$. We can 
then form the bundle $\xi:=P\times_{PU(\H)}Fred(\H)\to X$ and 
Atiyah and Segal define 
\[
K^{0}(X)_{f,AS}:=\pi_{0}\Gamma(\xi\to X),
\]
the group of homotopy classes of sections of $\xi\rightarrow X$.
In the appendix, we give details on how to use the symmetric 
spectrum model of $K$-theory due to Joachim \cite{joachim-ko} to 
obtain a map $K(\ZZ,3)\rightarrow BGL_1 K$ which is very much in 
the spirit of Atiyah and Segal.

The case of twisted real $K$-theory can be handled in the same way. 
As in the complex case, we are interested in those twists associated 
to a map $f$ representing a cohomology class $H^{2}(X;\Z/2)$. Thus 
given a map $f:X\to K(\ZZ/2,2)$ representing a cohomology class in 
$H^{2}(X;\Z/2)$ we can define the $f$-twisted real $K$-theory as 
$KO(X)_{i\circ f}$, where 
\[
i:K(\ZZ/2,2)\to BGL_{1}KO
\] 
is the inclusion map.

Therefore construction of 
twisted complex $K$-theory associated to integral cohomology 
classes in $\Hoh^3(X,\ZZ)$ amounts to a pointed map 
\[
j:K(\ZZ,3)\rightarrow BGL_1 K.
\]
Similarly, a construction of twisted real $K$-theory associated to 
cohomology classes in $\Hoh^2(X,\ZZ/2)$ 
amounts to a pointed map 
\[
j:K(\ZZ/2,2)\rightarrow BGL_1 KO.
\]
On the other hand, a construction for $r$-torsion integral 
classes in the complex case
is determined by a pointed map 
\[
j_r:K(\ZZ/r,2)\rightarrow BGL_1 K.
\]
A construction of all integral classes yields one for the 
$r$-torsion ones by composition with the Bockstein 
$\beta:K(\ZZ/r,2)\rightarrow K(\ZZ,3)$.

Our main theorem says there are no higher twists of complex 
$K$-theory on $K(\ZZ,3)$ or $K(\ZZ/r,2)$. In the real case we 
establish the nonexistence of higher twists of real 
$K$-theory on $K(\ZZ/2,2)$. Let $p:BGL_1 K\rightarrow K(\ZZ,3)$  and 
$q:BGL_{1}KO\to K(\ZZ,2)$ be the projection maps.

\begin{theorem}\label{uniqueness}
    The map $p$ induces isomorphisms
    \begin{align*}
\sigma: [K(\ZZ,3),BGL_1 K]& \rightarrow [K(\ZZ,3),K(\ZZ,3)]\cong \ZZ \\
 \sigma_r:[K(\ZZ/r,2),BGL_1 K]&\rightarrow[K(\ZZ/r,2),K(\ZZ,3)]\cong\ZZ/r.
\end{align*}   
The map $q$ induces an isomorphism
    \[
    \tau:[K(\ZZ/2,2),BGL_1 KO]\rightarrow[K(\ZZ/2,2),K(\ZZ/2,2)]\cong \ZZ/2.
    \]
    \begin{proof}
        The inclusion of the $K(\ZZ,3)$ component into $BGL_1 K$ gives 
        surjectivity. So, it suffices to prove that $\sigma$ is injective. 
        In other words, we wish to show that if $\sigma(j)=0$, then $j$ 
        is null-homotopic. But, if $\sigma(j)=0$, then the map
        \[
        K(\ZZ,3)\rightarrow K(\ZZ/2,1)\times K(\ZZ,3)\times BBSU_{\otimes}
        \]
        is homotopically trivial on the first two components. 
        By Proposition~\ref{prop:homtriv}, it is also trivial on the third 
        component. The statement for $\sigma_r$ is similar. The statement 
        for $\tau$ is proved in the same way by using 
        Proposition~\ref{prop:homtriv real}.
    \end{proof}
\end{theorem}

The previous theorem shows that in particular, any definition of complex 
twisted $K$-theory arising through a map 
\[
K(\ZZ,3)\to BGL_1 K
\]
agrees, up to multiplication of an integer, with the definition given 
above. For a given definition this integer can be obtained by 
determining the differential $d_{3}$ in the Atiyah-Hirzebruch spectral 
sequence computing the twisted equivariant $K$-groups. Equivalently, 
we can determine this integer by computing the twisted $K$-groups on the 
sphere $S^{3}$ for a generator $\alpha\in\Hoh^3(S^3,\ZZ)\cong \ZZ$. A 
similar situation occurs for twisted $K$-theory associated to $r$-torsion 
integral classes. For the case of twisted real $K$-theory the situation 
simplifies. In this case, by Theorem \ref{uniqueness} we have 
$[K(\ZZ/2,2),BGL_1 KO]\cong \ZZ/2$. In particular, 
any two non-trivial definitions of twisted real $K$-theory arising 
through a map 
\[
K(\ZZ/2,2)\to BGL_{1}KO
\]
must coincide.

\begin{corollary}
    Let $j:K(\ZZ,3)\rightarrow BGL_1 K$ be a pointed map. Then,
    \begin{equation*}
        \Omega j:K(\ZZ,2)\rightarrow GL_1 K
    \end{equation*}
    sends a line bundle $L$ on $X$ to the equivalence 
    $K(X)\stackrel{\sim}{\rightarrow} K(X)$ given by tensoring 
    with $L^{\otimes \sigma(j)}$. Similarly, if 
    $j:K(\ZZ/2,2)\rightarrow BGL_1 KO$ is a pointed map, 
    then $\Omega j$ is the auto-equivalence of $KO$
    given by tensoring with the $\tau(j)$-th power of real 
    line bundles.
    \begin{proof}
        This is true by construction for the morphism 
        $K(\ZZ,3)\rightarrow BGL_1 K$, which has $\sigma(j)=1$, 
        constructed in \cite{abg}. Thus, it is true for all other 
        definitions.
    \end{proof}
\end{corollary}

\section{Appendix: a geometric model of twisted $K$-theory}\label{appendix}

We outline how the symmetric spectrum model of $K$-theory (or $KO$-theory) 
due to Joachim \cite{joachim-ko} may be used to twist $K$-theory in a fashion 
that is nice from both the geometric and homotopical perspectives. The original
model for twisted $K$-theory spectra is due to Atiyah and Segal 
\cite{atiyah-segal}, but it is not easy to see directly how it fits into 
the homotopical framework of twists in the sense of \cite{abghr}. On the other 
hand, in \cite{abghr} and \cite{abg} it is hard to see the concrete analysis 
and geometry which were the original foundation for twisted $K$-theory.
Joachim's spectrum provides a vantage where both views may be appreciated.

In his paper, Joachim works with real periodic $K$-theory, but no 
alterations except replacing $8$ by $2$ at various places
are required to apply the same arguments to complex periodic 
$K$-theory. For simplicity, we present the complex case.

Let $\H$ be a fixed infinite dimensional separable Hilbert space, and 
let $\H_*=\H_0\oplus\H_1$ be the $\ZZ/2$-graded Hilbert space
with $\H_0=\H_1=\H$. The space $U(\H)$ is the group of all unitary 
operators on $\H$, equipped with the norm topology;
it is a contractible space by Kuiper's theorem.
The quotient of $U(\H)$ by the subgroup $U(1)$ of diagonal operators is 
$PU(\H)\we K(\ZZ,3)$. If $F:\H\rightarrow\H$ is an operator and 
$P\in PU(\H)$, then $P^{-1}FP$ is another operator.
Let $\F^1(\H_*)$ be the space of self-adjoint odd Fredholm 
operators on $\H_*$.
Thus, any element of $F$ of $\F^1(\H_*)$ can be represented by 
a matrix
\begin{equation*}
    \begin{pmatrix}
        0   &   \tilde{F}\\
        \tilde{F}^* &   0
    \end{pmatrix},
\end{equation*}
where $\tilde{F}$ is a Fredholm operator on $\H$. The group $PU(\H)$ 
acts continuously on $\F^1(\H_*)$ by
\begin{equation*}
    P\cdot    \begin{pmatrix}
        0   &   \tilde{F}\\
        \tilde{F}^* &   0
    \end{pmatrix}=    \begin{pmatrix}
        0   &   P^{-1}\tilde{F}P\\
        P^{-1}\tilde{F}^*P &   0
    \end{pmatrix},
\end{equation*}
where we purposely confuse $P$ with any operator in $U(\H)$ 
representing $P$.

If one is only interested in twisted $K$-groups, then this is 
already enough setup to do so. If $P$ is a principal 
$PU(\H)_{\mathrm{norm}}$-bundle on $X$, Atiyah and Segal define 
$K^0_P(X)$ as the group of homotopy classes
of sections of the associated bundle
\begin{equation*}
    P\times_{PU(\H)}\F^1(\H_*)\rightarrow X
\end{equation*}
with fiber $\F^1(\H_*)$. However, we are of course interested in 
an entire spectrum.

To describe twisted $K$-theory spectrum of Joachim, we recall the 
Clifford algebra $Cl(n)$, the complex Clifford algebra of $\CC^n$ 
equipped with the quadratic form
\begin{equation*}
    q((z_1,\ldots,z_n))=-\sum_{i=1}^n z_i^2.
\end{equation*}
There are canonical isomorphisms
\begin{equation*}
    Cl(p)\widehat{\otimes}Cl(q)\rightarrow Cl(p+q)
\end{equation*}
for $p,q\geq 1$, where $\widehat{\otimes}$ denotes the 
$\ZZ/2$-graded tensor product. For details on Clifford algebras, 
consult \cite{lawson-michelsohn}.

If $\J_*$ is a $\ZZ/2$-graded Hilbert space module for $Cl(n)$, 
we let $\F^1_{Cl(n)}(\J_*)$ denote the space of odd self-adjoint 
Fredholm operators on $\J_*$ which are $Cl(n)$-module morphisms when 
$n$ is even or the complement of the two contractible components 
of this space identified in \cite{atiyah-singer} when $n$ is odd.

Let
\begin{equation*}
    \H(n)_*=(Cl(1)\widehat{\otimes}\H_*)^{\hat{\otimes}n}.
\end{equation*}
Then, $\H(n)_*$ is naturally a graded $Cl(n)$-module. Joachim 
shows that the multiplication maps
\begin{gather*}
    \mu_{p,q}:\F^1_{Cl(p)}(\H(p)_*)\times\F^1_{Cl(q)}
    (\H(q)_*)\rightarrow\F^1_{Cl(p+q)}(\H(p+q)_*),\\
    (F,G)\mapsto F\star G=F\widehat{\otimes}Id+Id\widehat{\otimes}G
\end{gather*}
are continuous. The maps $\mu_{p,q}$ are 
$\Sigma_p\times\Sigma_q$-equivariant, where $\Sigma_n$ acts
naturally on $\F^1_{Cl(n)}(\H(n)_*)$.

To create based maps, let $K_n=\F^1_{Cl(n)}(\H(n))_+$, topologized 
as in \cite[Section~3]{joachim-ko} so that 
$\F^1_{Cl(n)}(\H(n))\rightarrow K_n$ is a homotopy equivalence. Then, 
the $\mu_{p,q}$ induce continuous maps 
$K_p\wedge K_q\rightarrow K_{p+q}$ for $p,q\geq 1$.

Let $PU(\H)$ act on $\F^1_{Cl(n)}(\H(n)_*)$ in the natural way, 
through the diagonal action of $U(\H)$ on $\H_*^{\widehat{\otimes} n}$. 
This extends to a continuous action of $PU(\H)$ on $K_n$, 
where $PU(\H)$ fixes the basepoint.

\begin{proposition}
    The natural $PU(\H)$ actions on the spaces of Joachim's 
    spectrum $K$ make $K$ into a $PU(\H)$-spectrum in the sense 
    that the actions are compatible with the multiplication maps and 
    the symmetric group actions.
    \begin{proof}
        This is clear as $PU(\H)$ acts diagonally on 
        $\H_*^{\widehat{\otimes}n}$.
    \end{proof}
\end{proposition}

\begin{corollary}\label{cor action PU(H)}
    The conjugation action of $PU(\H)$ on Joachim's spectrum $K$ determines an 
    $A_{\infty}$-map $PU(\H)\rightarrow GL_1 K$ which deloops to
    a map $K(\ZZ,3)\we BPU(\H)\rightarrow BGL_1K$.
\end{corollary}


\begin{bibdiv}
  \begin{biblist}


\bib{Adamsbook}{book}{
author={Adams, J. F.},
title={Stable homotopy and generalised homology},
note={Chicago Lectures in Mathematics},
publisher={University of Chicago Press},
place={Chicago, Ill.},
date={1974},
pages={x+373},
}

\bib{Adams}{article}{
author={Adams, J. F.},
author={Priddy, S. B.},
title={Uniqueness of $B{\rm SO}$},
journal={Math. Proc. Cambridge Philos. Soc.},
volume={80},
date={1976},
number={3},
pages={475--509},
issn={0305-0041},
}

\bib{anderson-hodgkin}{article}{
author={Anderson, D. W.},
author={Hodgkin, Luke},
title={The $K$-theory of Eilenberg-MacLane complexes},
journal={Topology},
volume={7},
date={1968},
pages={317--329},
issn={0040-9383},
}

\bib{abg}{article}{
author={Ando, Matthew},
author={Blumberg, Andrew J.},
author={Gepner, David},
title={Twists of $K$-theory and TMF},
conference={
title={Superstrings, geometry, topology, and $C^\ast$-algebras},
},
book={
series={Proc. Sympos. Pure Math.},
volume={81},
publisher={Amer. Math. Soc.},
place={Providence, RI},
},
date={2010},
pages={27--63},
}

\bib{abghr}{article}{
author = {Ando, M.},
author = {Blumberg, A. J.},
author = {Gepner, D. J.},
author = {Hopkins, M. J.},
author = {Rezk, C.},
title = {Units of ring spectra and Thom spectra},
journal = {ArXiv e-prints},
eprint = {http://arxiv.org/abs/0810.4535}, 
}

\bib{antieau-williams}{article}{
author = {Antieau, B.},
author = {Williams, B.},
title = {The period-index problem for twisted topological K-theory},
journal = {ArXiv e-prints},
eprint = {http://arxiv.org/abs/1104.4654},
}


\bib{atiyah-segal}{article}{
author={Atiyah, Michael},
author={Segal, Graeme},
title={Twisted $K$-theory},
journal={Ukr. Mat. Visn.},
volume={1},
date={2004},
number={3},
pages={287--330},
issn={1810-3200},
}

\bib{atiyah_twisted_2006}{article}{
author={Atiyah, Michael},
author={Segal, Graeme},
title={Twisted $K$-theory and cohomology},
conference={
title={Inspired by S. S. Chern},
},
book={
series={Nankai Tracts Math.},
volume={11},
publisher={World Sci. Publ., Hackensack, NJ},
},
date={2006},
pages={5--43},
}

\bib{atiyah-singer}{article}{
author={Atiyah, M. F.},
author={Singer, I. M.},
title={Index theory for skew-adjoint Fredholm operators},
journal={Inst. Hautes \'Etudes Sci. Publ. Math.},
number={37},
date={1969},
pages={5--26},
issn={0073-8301},
}

\bib{bcmms}{article}{
author={Bouwknegt, Peter},
author={Carey, Alan L.},
author={Mathai, Varghese},
author={Murray, Michael K.},
author={Stevenson, Danny},
title={Twisted $K$-theory and $K$-theory of bundle gerbes},
journal={Comm. Math. Phys.},
volume={228},
date={2002},
number={1},
pages={17--45},
issn={0010-3616},
}

\bib{cartan-determination}{article}{
    AUTHOR = {Cartan, H.},
    TITLE = {D{\'e}termination des alg{\`e}bres {$\Hoh_*(\pi,n;\ZZ)$}},
    JOURNAL = {S{\'e}minaire H. Cartan},
    VOLUME = {7},
    NUMBER = {1},
    PAGES = {11-01--11-24},
    PUBLISHER = {Secr{\'e}tariat math{\'e}matique},
    ADDRESS = {Paris},
    YEAR = {1954/1955},
}


\bib{donovan_karoubi}{article}{
author={Donovan, P.},
author={Karoubi, M.},
title={Graded Brauer groups and $K$-theory with local coefficients},
journal={Inst. Hautes \'Etudes Sci. Publ. Math.},
number={38},
date={1970},
pages={5--25},
issn={0073-8301},
}

\bib{fht}{article}{
author={Freed, Daniel S.},
author={Hopkins, Michael J.},
author={Teleman, Constantin},
title={Twisted equivariant $K$-theory with complex coefficients},
journal={J. Topol.},
volume={1},
date={2008},
number={1},
pages={16--44},
issn={1753-8416},
}

\bib{gomez}{article}{
author={G{\'o}mez, Jos{\'e} Manuel},
title={On the nonexistence of higher twistings over a point for Borel
cohomology of $K$-theory},
journal={J. Geom. Phys.},
volume={60},
date={2010},
number={4},
pages={678--685},
issn={0393-0440},
}

\bib{griffiths-morgan}{book}{
author={Griffiths, Phillip A.},
author={Morgan, John W.},
title={Rational homotopy theory and differential forms},
series={Progress in Mathematics},
volume={16},
publisher={Birkh\"auser Boston},
place={Mass.},
date={1981},
pages={xi+242},
isbn={3-7643-3041-4},
}

\bib{joachim-ko}{article}{
author={Joachim, Michael},
title={A symmetric ring spectrum representing $KO$-theory},
journal={Topology},
volume={40},
date={2001},
number={2},
pages={299--308},
issn={0040-9383},
}


\bib{karoubi}{article}{
author={Karoubi, Max},
title={Twisted $K$-theory---old and new},
conference={
title={$K$-theory and noncommutative geometry},
},
book={
series={EMS Ser. Congr. Rep.},
publisher={Eur. Math. Soc., Z\"urich},
},
date={2008},
pages={117--149},
}


\bib{lawson-michelsohn}{book}{
author={Lawson, H. Blaine, Jr.},
author={Michelsohn, Marie-Louise},
title={Spin geometry},
series={Princeton Mathematical Series},
volume={38},
publisher={Princeton University Press},
place={Princeton, NJ},
date={1989},
pages={xii+427},
isbn={0-691-08542-0},
}

\bib{htt}{book}{
author={Lurie, Jacob},
title={Higher topos theory},
series={Annals of Mathematics Studies},
volume={170},
publisher={Princeton University Press},
place={Princeton, NJ},
date={2009},
pages={xviii+925},
isbn={978-0-691-14049-0},
isbn={0-691-14049-9},
}

\bib{mst}{article}{
author={Madsen, I.},
author={Snaith, V.},
author={Tornehave, J.},
title={Infinite loop maps in geometric topology},
journal={Math. Proc. Cambridge Philos. Soc.},
volume={81},
date={1977},
number={3},
pages={399--430},
issn={0305-0041},
}

\bib{may-sigurdsson}{book}{
author={May, J. P.},
author={Sigurdsson, J.},
title={Parametrized homotopy theory},
series={Mathematical Surveys and Monographs},
volume={132},
publisher={American Mathematical Society},
place={Providence, RI},
date={2006},
pages={x+441},
isbn={978-0-8218-3922-5},
isbn={0-8218-3922-5},
}

\bib{McGibbon}{article}{
author={McGibbon, C. A.},
author={M{\o}ller, J. M.},
title={On spaces with the same $n$-type for all $n$},
journal={Topology},
volume={31},
date={1992},
number={1},
pages={177--201},
issn={0040-9383},
}

\bib{meier}{article}{
author={Meier, W.},
title={Complex and real $K$-theory and localization},
journal={J. Pure Appl. Algebra},
volume={14},
date={1979},
number={1},
pages={59--71},
issn={0022-4049},
}

\bib{milnor}{article}{
author={Milnor, J.},
title={On axiomatic homology theory},
journal={Pacific J. Math.},
volume={12},
date={1962},
pages={337--341},
issn={0030-8730},
}

\bib{rosenberg}{article}{
author={Rosenberg, Jonathan},
title={Continuous-trace algebras from the bundle theoretic point of view},
journal={J. Austral. Math. Soc. Ser. A},
volume={47},
date={1989},
number={3},
pages={368--381},
issn={0263-6115},
}

\bib{wang}{article}{
author={Wang, Bai-Ling},
title={Geometric cycles, index theory and twisted $K$-homology},
journal={J. Noncommut. Geom.},
volume={2},
date={2008},
number={4},
pages={497--552},
issn={1661-6952},
}

\bib{weibel}{book}{
author={Weibel, Charles A.},
title={An introduction to homological algebra},
series={Cambridge Studies in Advanced Mathematics},
volume={38},
publisher={Cambridge University Press},
place={Cambridge},
date={1994},
pages={xiv+450},
isbn={0-521-43500-5},
isbn={0-521-55987-1},
}

\bib{yosimura}{article}{
author={Yosimura, Zen-ichi},
title={A note on complex $K$-theory of infinite CW-complexes},
journal={J. Math. Soc. Japan},
volume={26},
date={1974},
pages={289--295},
issn={0025-5645},
}


  \end{biblist}
\end{bibdiv}
\end{document}